\documentclass[ijoc,nonblindrev]{informs3} 

\OneAndAHalfSpacedXII 

\usepackage{natbib}
 \bibpunct[, ]{(}{)}{,}{a}{}{,}%
 \def\bibfont{\small}%
\TheoremsNumberedThrough     

\EquationsNumberedThrough    
\usepackage[T1]{fontenc} 
\usepackage{lmodern} 
\usepackage{amsmath,amsfonts,amssymb}

\usepackage{xcolor}
\usepackage{graphicx}
\usepackage{multirow}
\usepackage{booktabs}
\usepackage{subcaption}
\usepackage{threeparttable}
\usepackage{authblk}
\usepackage[colorlinks=true, allcolors=blue]{hyperref}
\usepackage{color}
\usepackage[colorinlistoftodos,prependcaption,textsize=tiny]{todonotes}

\newcommand{\init}{\texttt{root}}
\newcommand{\fin}{\texttt{end}}

\DeclareMathOperator{\thepace}{pace}
\newcommand{\OGpace}[1][]{\thepace^{OG}_{#1}}
\newcommand{\LBpace}[1][]{\thepace^{LB}_{#1}}

\DeclareMathOperator{\npace}{npace}

\newcommand{\NLBpace}[1][]{\npace^{LB}_{#1}}

\begin{document}


\RUNAUTHOR{Ghaddar et al.}

\RUNTITLE{Learning for Spatial Branching}

 \TITLE{Learning for Spatial Branching: An Algorithm Selection Approach}

\ARTICLEAUTHORS{
\AUTHOR{Bissan Ghaddar}{\AFF{Corresponding author. Ivey Business School, Western University, London, Ontario, Canada, \EMAIL{bghaddar@ivey.ca}}}
\AUTHOR{Ignacio G\'omez-Casares, Julio Gonz\'alez-D\'iaz, Brais Gonz\'alez-Rodr\'iguez, Beatriz Pateiro-L\'opez}{\AFF{CITMAga (Galician Center for Mathematical Research and Technology). Department of Statistics, Mathematical Analysis and Optimization and MODESTYA Research Group, University of Santiago de Compostela.}}
\AUTHOR{Sof\'ia Rodr\'iguez-Ballesteros}{\AFF{CITMAga (Galician Center for Mathematical Research and Technology) }}
}

\ABSTRACT{The use of machine learning techniques to improve the performance of branch-and-bound optimization algorithms is a very active area in the context of mixed integer linear problems, but little has been done for non-linear optimization. To bridge this gap, we develop a learning framework for spatial branching and show its efficacy in the context of the Reformulation-Linearization Technique for polynomial optimization problems. The proposed learning is performed offline, based on instance-specific features and with no computational overhead when solving new instances. Novel graph-based features are introduced, which turn out to play an important role for the learning. Experiments on different benchmark instances from the literature show that the learning-based branching rule significantly outperforms the standard rules.}

\maketitle

\section{Introduction} 


In recent years, there has been a substantial rise on the use of learning techniques to improve the performance of branch-and-bound schemes. Most of the focus so far has been set on mixed integer linear programming (MILP) problems and, in particular, on learning to choose the branching variable (see \cite{Lodi2017} and \cite{Bengio2021} and references therein). Since branch-and-bound schemes are also at the core of most global optimization algorithms for solving non-linear and non-convex optimization problems, it is natural to wonder to what extent the insights gained for (discrete) MILP problems can be successfully transferred to the (continuous) non-linear setting. To this end, in this paper we study variable selection for spatial branching in non-linear programming (NLP) problems.

We develop our analysis for polynomial optimization (PO) problems, a subclass of NLP problems known to be $\mathcal{NP}$-hard in the general case. PO has a wide range of practical applications in the context of control, energy and water networks, process engineering, facility location, economics and equilibrium, and finance. PO generalizes several special cases that have been thoroughly studied in optimization, including mixed binary linear optimization, convex/non-convex quadratic optimization, and complementarity problems. 

A key factor of the efficiency of state-of-the-art global solvers in non-linear optimization is the design of the branch-an-bound algorithm and, in particular, the branching rules. Inspired by the recent literature on ``learning to branch'', we explore ways to harness learning-based approaches to improve the performance of the branch-and-bound search for PO and, more specifically, we focus on variable selection in spatial branching. Among the various approaches that have been proposed for solving problems in PO, here we build upon the branch-and-bound scheme embedded in the Reformulation-Linearization Technique (RLT) introduced by \cite{Sherali1992}. Although we have framed our contribution in the context of the RLT technique for PO, the insights are extensible to any branch-and-bound technique designed to find global optima for any class of NLP problems.


\section{Contribution to the Literature} 

Machine learning (ML) and the neighbouring field of statistical learning have both been effectively leveraged in several areas to develop algorithms that learn from observing the performance of different tasks. While mathematical optimization intrinsically lies at the core of learning methods, recent years have seen a rise of research in the opposite direction: learning and predicting the best optimization approaches in different settings. In particular, a very active area of research concerns the use of ML techniques in the algorithmic design of optimization solvers, independently of the application domain in which such algorithms are used. A substantial part of the associated literature, to which this paper belongs, focuses on improving the performance of branch-and-bound schemes.

Because of the wide ranging applications and the ease of use of state-of-the-art solvers, mixed integer linear programming represents an excellent testing ground for the aforementioned research. Further, since these problems are solved using branch-and-bound techniques, the richness of the data generated in the process naturally lends itself to learning. As already mentioned, variable selection in branch-and-bound schemes is probably the topic that has attracted the most attention so far and its outstanding importance for the performance of such schemes has long been recognized. 




The integration of learning in the algorithmic design of the branching decisions can involve different paradigms and approaches, and occur on various levels. Strong branching has proven to be the variable selection rule that leads to smallest trees \citep{Linderoth1999}, but the associated computational costs make it's overall performance not competitive with more sophisticated rules such as reliability branching \citep{Achterberg2005}. This motivates the research in \cite{Khalil2016} and \cite{Alvarez2017}, where the authors devise learning approaches to approximate strong branching without the associated computational overhead, and test their approaches on the MIPLIB instances.
Another natural learning approach related to branching is taken in \cite{Liberto2016} where the authors, inspired by the portfolio approaches to algorithm selection \citep{Gomes2001,Leyton2003}, learn a clustering-based classifier to switch between branching rules and again test their approach on MIPLIB instances. Because of the richness of the different learning paradigms, novel approaches in this area are constantly appearing such as the use of neural networks in \cite{Gasse2019} and \cite{Cappart2021}. They approximate strong branching by encoding the branching policies into a graph convolutional neural network which allows exploiting the natural bipartite graph representation of MILP problems, thereby reducing the amount of manual feature engineering. Then, they test their approach on four benchmarks of specific combinatorial optimization problems.

Although some work has been done on using learning for NLP, this area is still not as mature as the one for MILP. For instance, for non-convex quadratic programming problems with box constraints, \cite{Baltean2019} try to learn effective linear outer-approximations of semidefinite constraints. To generate the corresponding cuts, a neural network is used to select the most promising submatrices without the computational burden of solving semidefinite programming problems. For convex quadratic problems, \cite{Bonami2018} use ML to decide between applying branch and bound on the quadratic problem and applying it on an equivalent MILP reformulation.

Despite the extensive work done on learning to branch in MILP problems, it is an issue not yet well understood and a very active area of research. The understanding of variable selection is even more limited in spatial branching for NLPs, where there has not been much research in the topic (an exception is \cite{Belotti2009}), much less in the learning context. 

Our research aims to bridge this gap by proposing a learning framework for spatial branching. We start by noticing that even relatively similar branching rules exhibit a lot of variability in their relative performance to one another across instances of different test sets. This is not something specific to this setting since, as argued in \cite{Liberto2016} and references therein, work in portfolio algorithms has shown that there is often no single approach that performs best on all instances. Thus, it is important to know to what extent one can learn which rule will perform best on a given instance and this is the main objective of this paper. We follow an algorithm selection approach in the spirit of the one used in \cite{Liberto2016} for MILP problems, but our learning takes place in a regression framework instead of a clustering one. In the present work, learning is done offline and the branching rule for a new instance is selected based on its underlying features. Given the promising results we obtain, the proposed approach can be seen as a first step towards more sophisticated (possibly online) learning schemes for variable selection in NLP problems. Given that the learning takes place offline, our approach entails no computational overhead once the solving process starts, so there is no node-efficiency \emph{vs.} time-efficiency trade-off. We develop our analysis for PO problems by building upon the global solver RAPOSa \citep{Gonzalez-rodriguez2020}. We perform our numerical experiments using a comprehensive set of instances that builds upon well established benchmarks including randomly generated instances from \cite{Dalkiran2016}, MINLPLib instances \citep{minlplib}, and QPLIB instances \citep{qplib}, which illustrate that our approach can be successfully applied to general PO problems without focusing on a particular class of problems or structure.

The contribution of this paper is threefold: (1) developing a framework for learning to select the most efficient branching rule in non-linear optimization problems, (2) defining graph-based features and branching rules that significantly contribute to the outcome of the learning process, and (3) introducing a new performance measure which, by combining information related to both running time and optimality gap, allows for a smooth comparison in sets of instances of varying structure and difficulty.


The reminder of this paper is organized as follows. In Section~\ref{sec:motivation}, we present the core elements on which our learning framework stands: i) the RLT-based spatial branch-and-bound scheme for PO and ii) the set of branching rules for algorithm selection. In Section~\ref{sec:methodology}, we describe the learning methodology, motivate and define the novel key performance indicator (KPI), and describe the features used for learning. In Section~\ref{sec:results}, we show the performance of the learning approach. Finally,  conclusions and future research directions are drawn in Section~\ref{sec:conclusions}.
\section{Motivation and Problem Description}\label{sec:motivation} 
As we have discussed, the main goal of this paper is to provide a first study of the impact that learning may have in the context of spatial branching. In this section we describe the polynomial optimization setting in which we develop our analysis, some graphs and graph-related features for polynomial optimization problems, and the branching rules that constitute our portfolio for algorithm selection.

\subsection{Polynomial Optimization and RLT Technique} \label{sec:PO}
The Reformulation-Linearization Technique was originally developed in \cite{Sherali1992}. It was designed to find global optima in polynomial optimization problems of the following form:
\begin{equation}
\begin{split}
\text{minimize} & \quad \phi_0(\mathbf{x})\\
\text{subject to}  & \quad \phi_r(\mathbf{x})\geq \beta_r, \quad r=1,2,\ldots, R_1 \\
& \quad \phi_r(\mathbf{x})=\beta_r, \quad r=R_1+1,\ldots,R\\
& \quad \mathbf{x}\in\Omega \subset \mathbb{R}^n\text{,}
\end{split}
\tag{\textbf{PO}}
\label{eq:PO}
\end{equation}
where $N = \lbrace 1, \dots, n \rbrace$ denotes the set of variables, each $\phi_r(\mathbf{x})$ is a polynomial of degree $\delta_r \in \mathbb{N}$ and $\Omega = \lbrace \mathbf{x} \in \mathbb{R}^n: 0 \leq l_j \leq x_j \leq u_j < \infty, \, \forall j \in N \rbrace \subset \mathbb{R}^n$ is a hyperrectangle containing the feasible region. Then, $\delta=\max_{r \in \{0,\ldots,R\}} \delta_r$ is the degree of the problem and $(N, \delta)$ represents all possible monomials of degree $\delta$.

RLT is based on the construction of a linear relaxation of the polynomial problem, which is then solved using a branch-and-bound scheme. This is accomplished by replacing each monomial of the problem with a corresponding RLT variable. For example, associated to a monomial of the form $x_1x_2x_4$ one would define the RLT variable $X_{124}$. More generally, RLT variables are defined as
\begin{equation}
  X_J = \prod_{j \in J}x_j, 
  \label{eq:RLTidentity}
\end{equation}
where $J$ is a multiset containing the information about the multiplicity of each variable in the underlying monomial. Then, at each node of the branch-and-bound tree, one would solve the corresponding linear relaxation. In order to get a tighter relaxation and ensure convergence, new constraints, called bound-factor constraints, are also added. At each node, branching is performed on one of the variables that present some violation of the RLT-defining identities~\eqref{eq:RLTidentity}. Therefore, one of the crucial steps in any implementation of the RLT is the branching variable selection, which we discuss in Section~\ref{sec:branchingcriteria}.

\subsection{Graphs Associated to Polynomial Optimization Problems}\label{sec:graph}

The relevance of graph representations of general optimization problems has long been recognized and, as discussed in \cite{Bienstock2018} and \cite{Faenza2022}, polynomial optimization problems are not an exception. Our interest on these graphs is twofold: first, as a means to define new branching rules and, second, to define features that can be taken as input by the learning process.

To each problem of the form \ref{eq:PO}, one can associate different graphs that capture  some characteristics of its structure. We introduce two different graphs, both based on the concept of intersection graphs as defined in \cite{Fulkerson1965}. The first one is the \textit{variables intersection graph}, VIG, analogous to the one used in \cite{Bienstock2018}. It is built by assigning a vertex to each variable and connecting two vertices if the corresponding variables appear together in some monomial of~\ref{eq:PO}. The second one is the \textit{constraints-monomials intersection graph}, CMIG, where we have one vertex for each monomial, one vertex for each constraint, and one for the objective function. We then connect one vertex associated to a monomial with one associated to a constraint (or the objective function) if the monomial appears in the constraint (or objective function).

Once the graphs are defined, we are interested in parameters that summarize some properties of their underlying structure. Below we describe several standard parameters that we use in different parts of our analysis (for further details refer, for instance, to \cite{latora2017}).

\begin{itemize}
 \item \textbf{Density or edge density}. Given a graph $G$, the edge density of $G$ is defined as the number of edges in $G$ divided by the number of edges in the complete graph with the same number of nodes.
 \item \textbf{Treewidth}. Given a graph $G$, the treewidth is the smallest width of any tree-decomposition of $G$ \citep{robertson1984}. It gives a measure of how close a graph is to a tree. Trees have treewidth~1.
 \item \textbf{Modularity}. Given a division of a graph $G$ into communities, the modularity of the division is a measure of the strength of that division (\emph{i.e.}, graphs with high modularity have dense connections between the nodes within communities but sparse connections between nodes in different communities). The modularity of $G$ is the maximum modularity over all the possible divisions into communities of the graph. Since the exact computation of this parameter is a hard computational problem, we use the greedy algorithm in \cite{Clauset2004} to compute an approximate value.
 \item \textbf{Transitivity}. The transitivity of a graph $G$ is the ratio between the number of triangles (subgraphs that are complete graphs of order 3) and the number of triads in $G$.
 \item \textbf{Eigenvector centrality or eigencentrality}. For each node in a graph, its centrality score is the corresponding value in the first eigenvector of the graph adjacency matrix.
\end{itemize}

We now move to the definition of different branching rules, two of which heavily rely on the eigencentralities of the nodes associated with the variables of the problem.

\subsection{Branching Rules}\label{sec:branchingcriteria}



State-of-the art solvers for non-linear optimization build upon two main themes to define their variable selection rules: i) adaptations of MILP methods such as strong branching and, most notably, reliability branching \citep{Achterberg2005}, as described in \cite{Belotti2009} and ii) contributions of the variables to the violations of non-linear terms using variants of the so-called violation transfer \citep{Tawarmalani2004, Belotti2009}. In the specific context of RLT, past literature has mainly focused on violations of the RLT-defining identities. In the seminal paper by \cite{Sherali1992}, starting from an optimal solution $(\bar X, \bar x)$ at a node in the branch-and-bound tree, the violation of each variable in~\ref{eq:PO} is computed according to the following expression:
\begin{equation}
\label{eq:maxcriteria}
 \theta_j = \max_{t \in \lbrace 1, \dots, \delta-1 \rbrace} \max_{J \in (N, \delta): \;\vert J \vert = t} \lbrace \vert \bar{X}_{J \cup \lbrace j \rbrace} - \bar{x}_j \bar{X}_{J} \vert \rbrace\text{.}
\end{equation}

Then, the variable that maximizes this value is chosen for branching. More recent contributions such as \cite{Sherali2012cuts,Sherali2012reduced} and \cite{Gonzalez-rodriguez2020} rely on variations of~\eqref{eq:maxcriteria} where the maxima are replaced with sums. In particular, \cite{Gonzalez-rodriguez2020} define various branching rules based on the formula:
\begin{equation}
\label{eq:sumweightscriteria}
  \theta_j = \sum_{t \in \lbrace 1, \dots, \delta-1 \rbrace} \sum_{J \in (N, \delta): \;\vert J \vert = t} w(j,J) \vert \bar{X}_{J \cup \lbrace j \rbrace} - \bar{x}_j \bar{X}_{J} \vert
\end{equation}
where $w(j,J)$ are the weights associated with the violations. In order to define our portfolio of branching rules we complement the original formulation in~\eqref{eq:maxcriteria} with various rules defined by using different weights inside the sum, which can be seen in Table~\ref{tab:rules}.\footnote{The reason behind choosing only one rule based on maxima is that previous analysis have shown that ``maxima''-based rules are normally outperformed by the ``sum''-based counterparts \citep{Gonzalez-rodriguez2020}.}

\begin{table}[!htbp]
\centering
\caption{Definition of the branching rules.}
\begin{tabular}{p{0.2\textwidth}|p{0.73\textwidth}}
\textbf{Branching rule} & \textbf{Definition}\\
\toprule\toprule
Maximum (\texttt{max}) & Eq.~\eqref{eq:maxcriteria}.\\
\hline
Sum (\texttt{sum}) & Eq.~\eqref{eq:sumweightscriteria}, with $w(j,J)=1$ for all $j$ and $J$.\\
\hline
Dual (\texttt{dual}) & Eq.~\eqref{eq:sumweightscriteria}, with $w(j,J)$ defined as the sum of the absolute values of the shadow prices of the problem constraints containing $J\cup \{j\}$.\\
\hline
Range (\texttt{range}) & Eq.~\eqref{eq:sumweightscriteria}, with $w(j,J)=\frac{\min{\lbrace\bar{u}_j - \bar{x}_j, \,\bar{x}_j - \bar{l}_j\rbrace}}{u_j-l_j}$, where $u_j$ and $l_j$ are the upper and lower bound of $x_j$ at the root node, \emph{i.e.} \ref{eq:PO}, and $\bar{u}_j$ and $\bar{l}_j$ are the bounds at the current node. \\
\hline
Eigencentrality \newline VIG (\texttt{eig-VI}) & Eq.~\eqref{eq:sumweightscriteria}, with $w(j,J)$ defined as the eigencentrality of $x_j$'s node in VIG.\\
\hline
Eigencentrality \newline CMIG (\texttt{eig-CMI}) & Eq.~\eqref{eq:sumweightscriteria}, with $w(j,J)$ defined as the eigencentrality of $x_j$'s node in CMIG ($0$ if $x_j$ never appears alone in a monomial).\\
\hline
\end{tabular}
\label{tab:rules}
\end{table}%


\begin{figure}[!htbp]
    \centering
    \begin{subfigure}{\textwidth}
		\centering
    \includegraphics[width=0.85\textwidth]{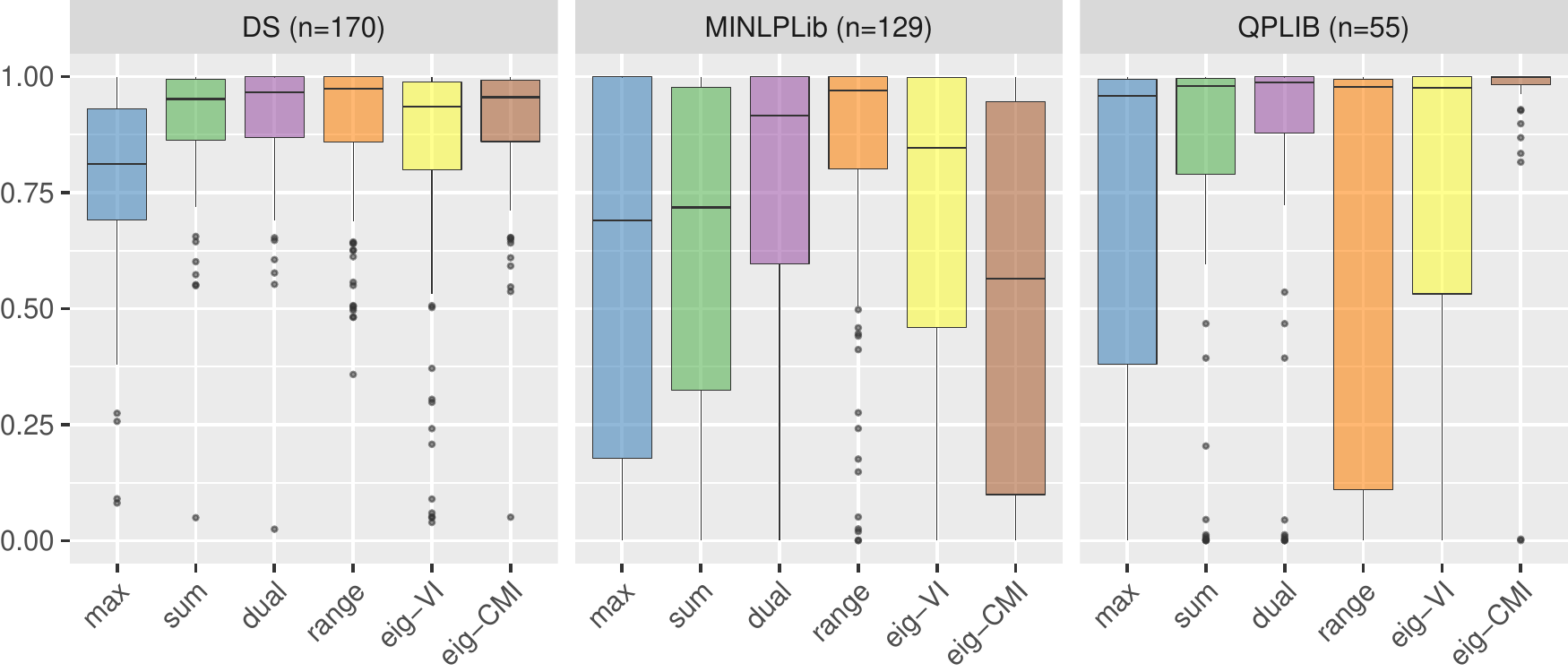}
    \caption{Boxplots of the normalized KPI for each of the branching rules.}
    \label{subfig:bp_comp}
    \end{subfigure}
    
    \medskip
    \begin{subfigure}{\textwidth}
		\centering
    \includegraphics[width=0.85\textwidth]{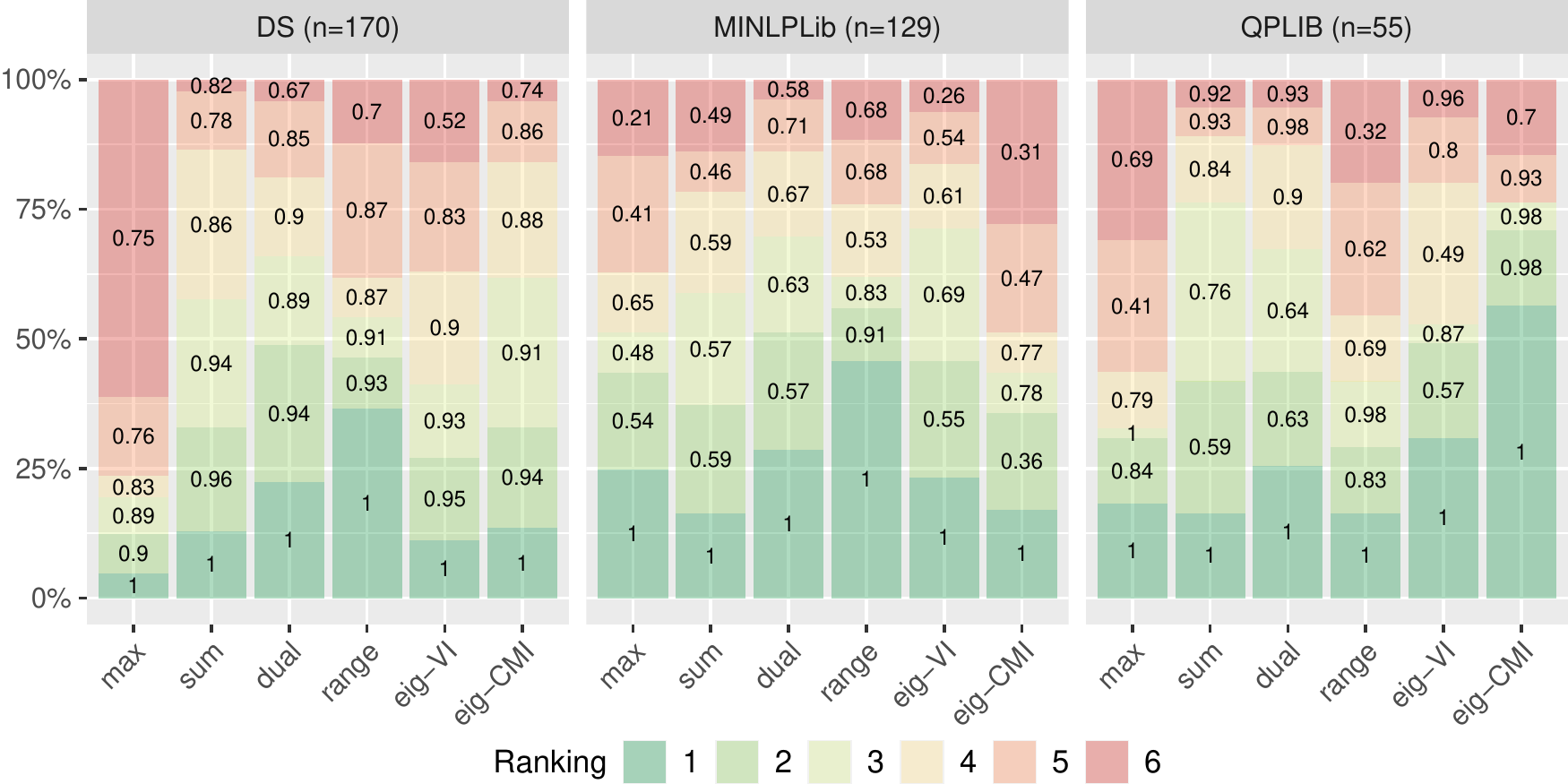}
    \caption{For each instance branching rules were ranked from 1 (best) to 6 (worst), according to the KPI. Stacked bar graphs represent percentage of problems in each rank position per branching rule. Numbers inside bars indicate the average value of the normalized KPI for the problems in each rank position.
}
    \label{subfig:bar_comp}
    \end{subfigure}
    \caption{Comparative analysis of the performance of the different branching rules.}
    \label{fig:branchingcriteria}
\end{figure}

Figure~\ref{fig:branchingcriteria} shows the results of a preliminary computational experiment to assess the performance of the different branching rules on PO problems. More precisely, we used three different sets of instances (refer to Section~\ref{sec:dataset} for details): DS~\citep{Dalkiran2016}, MINLPlib~\citep{minlplib}, and QPLIB~\citep{qplib}. For our comparison of the branching rules we use a key performance indicator (KPI) whose values are normalized to $[0, 1]$ (details are given in Section~\ref{sec:KPI}). The closer the normalized KPI is to~1 for a given branching rule, the better the performance of that rule.  Figure~\ref{subfig:bp_comp} shows, for the problems in each of the three sets of instances, the boxplots  of  this  normalized  KPI  using  each  of  the  described  branching  rules. As expected, no branching rule outperforms the rest in the different data sets. While the \texttt{sum} rule performs reasonably well in the DS instances, it does not for the ones in MINLPLib. In QPLIB, the \texttt{eig-CMI} rule clearly outperforms all other rules, but it would be a poor choice for many of the MINLP instances. These observations already hint at the promising potential of the integration of learning techniques in the branching rule selection process. In Figure~\ref{subfig:bar_comp} we rank, for each problem, the branching rules from 1 (best) to 6 (worst), according to the normalized KPI. As we observe from the learning process, it is important not only to choose the best possible branching rule, but also to ensure that, in case of mistake, the selected rule performs reasonably well. For instance, for those MINLPLib instances where the \texttt{max} rule is not the best one, its performance is on average significantly worse than other methods. In those cases it would be preferable to choose, for instance, the \texttt{range} rule, even when it is not the best possible choice. 



\section{Proposed Methodology}\label{sec:methodology}
\subsection{Overview of the Learning Framework}\label{sec:framework}
Our learning framework is related to the approach for switching heuristics in \cite{Liberto2016}, where the authors use a cluster-based method to learn to choose from a portfolio of branching rules in MILP problems. In this paper, we develop a different, regression-based, learning technique. The main elements of our setting, some of which are novel with respect to past literature, are the following:
\begin{itemize}
    \item \textbf{Spatial branching}. We separate from previous literature by studying branching variable selection in non-linear problems.
    \item \textbf{Algorithm portfolio}. The goal is learning to choose from the branching rules defined in Section~\ref{sec:branchingcriteria}.
    \item \textbf{Data collection}. Learning is performed on a rich set of very diverse instances (see Section~\ref{sec:dataset}). All instances are solved with all branching rules, gathering all the relevant data on their performance. Thus, all the learning takes place offline, with no computational overhead when solving new instances (beyond the computation of the features of the given instance).
    \item \textbf{Graph-based branching rules and features}. A novel part of our approach is the use of graphs to define both branching rules and features (see Tables~\ref{tab:rules} and~\ref{table:features}, respectively). Remarkably, both of them turn out to play an important role in the results.
    \item \textbf{Learning technique}. It consists of running individual regressions for each branching rule, with the goal of predicting its relative performance with respect to the remaining ones. Then, given a new instance, the rule with the best predicted relative performance is chosen.
    \item \textbf{KPI.} We present a novel KPI, which we refer to as pace, and that is discussed in detail in Section~\ref{sec:KPI}. The pace allows to jointly study the performance in instances solved by all methods (where running time is normally used) and instances not solved by at least one method (where optimality gap is normally used).
    \item \textbf{Learning on static features.} Learning uses instance-specific features, with no dynamic information being gathered on a node by node basis along the branch-and-bound tree. Thus, the resulting learned branching rule is static, in the sense that it selects the same branching rule in each and every node of the tree associated to a given instance. This is different from most of the literature in MILP problems where dynamic features are often included (see, for instance, \cite{Liberto2016}, \cite{Khalil2016}, and \cite{Alvarez2017}).
\end{itemize}

It is worth noting that, despite the variability in performance illustrated in Figure~\ref{fig:branchingcriteria}, all the rules in our portfolio revolve around the same baseline computation (violations of RLT-defining identities). This common core probably limits to what extent they can complement each other and, as a consequence, may result in learned rules with less potential to significantly outperform the original ones. Similarly, the choice to perform learning only with static features, may result in learning a less effective branching rule. Because of the two aforementioned points, our analysis provides some sort of lower bound on the ``potential for learning'' in spatial branching. One of our main contributions is to show that, even under these ``limitations'', there is much potential to improve the overall performance of the branch-and-bound algorithm by conducting offline learning on the branching rule.

\subsection{A Novel Key Performance Indicator}\label{sec:KPI} 
Carrying out a fair comparison of the performance of different algorithms on instances of varying difficulty often involves some challenges. Ideally, one would like to have a KPI that allows to compare the performance on all instances together, but the most widely used KPIs, running time and optimality gap, fail to do so. To illustrate this, let's define the following subsets of the set of instances: 
\begin{itemize}
	\item $A=\,$``Instances solved by all algorithms within the time limit'',
	\item $B=\,$``Instances solved by some but not all algorithms within the time limit'', and
	\item $C=\,$ ``Instances not solved by any algorithms within the time limit''. 
\end{itemize}
The running time is a good KPI in $A$ and uninformative in $C$, whereas the optimality gap is a good KPI in $C$ and uninformative in $A$. Further, neither the time nor the gap allows to differentiate between algorithms for all instances in $B$, since the running time does not distinguish between algorithms that have reached the time limit (but with different gaps) and the optimality gap does not distinguish between algorithms that have solved the given instance (but with different running time). Because of this, it is customary to split the performance analysis separating the problems in $A\cup B$ and $C$ or $A$ and $C\cup B$. For instance, in \cite{Alvarez2017} they ``separated the problems that were solved by all methods from the problems that were not solved by at least one of the compared methods''. An important limitation of such approaches is that the sets $A$, $B$ and $C$ depend on the algorithms to be compared, so the results may substantially change if a new algorithm is added to the study or an existing one is removed. In \cite{Gupta2020}, when discussing limitations of their work stemming from the $\mathcal{NP}$-hard  nature of MILP solving, the authors acknowledge that ``One can consider the primal-dual bound gap after a time limit as an evaluation metric for the bigger instances, but this is misaligned with the solving time objective''.

We propose a novel approach to define KPIs that allows to jointly compare all instances while being equivalent to the running time in $A$, mimicking the optimality gap in $C$, and being able to completely differentiate between algorithms in $B$. Thus, this new family of KPIs overcomes all the issues mentioned above. First, recall the standard definition of the optimality gap at the end of the algorithm, $OG^\fin$, as a function of the lower and upper bounds ($LB$ and $UB$, respectively):
\[
OG^\fin=\frac{UB^\fin-LB^\fin}{UB^\fin+\varepsilon},
\]
where $\varepsilon$ is a small constant needed to avoid divisions by zero. We can now define a KPI based on the pace at which the optimality gap is closed during the execution:
\[
\OGpace=\frac{\text{time}}{OG^\fin-OG^\init+\varepsilon}.
\]
Note that, in our setting, $OG^\init$ is common for all branching rules. So defined, $\OGpace$ represents the time needed to improve the optimality gap in one unit. Suppose now that we are comparing two algorithms according to their paces $\OGpace[1]$ and $\OGpace[2]$:
\begin{itemize}
    \item In set $A$, $OG^\fin$ is the same for all rules. Hence, $\OGpace[1]/\OGpace[2]=\text{time}_1/\text{time}_2$.
    \item In set $C$, the running time is the same for all rules. Hence, $\OGpace[1]/\OGpace[2]$ is of the form $(OG^\fin_2+K)/(OG^\fin_1+K)$, where $K$ is a constant term.
    \item In set $B$, $\OGpace$ may recognize the difference between two rules that have solved an instance if they have different running times and, similarly, between two rules that have not solved an instance if they have different gaps. Further, as desired, rules that have solved an instance will have a better pace than those that have not.
\end{itemize}
Given the above properties, $\OGpace$ stands out as natural candidate for a streamlined comparison of branching rules, jointly for all instances. Unfortunately, although branch-and-bound schemes always compute a lower bound at the root node, an initial upper bound is often not available. Yet, since the main goal of the branching rules is often to deliver a faster increase of the lower bounds, it is natural to define a variant of $\OGpace$ that is based on the pace at which the lower bound increases:
\begin{equation}
    \LBpace=\frac{\text{time}}{LB^\fin-LB^\init+\varepsilon}.
\label{eq:epsilon}
\end{equation}

Note that $\LBpace$ is always well defined and measures the amount of time needed to improve the lower bound in one unit. This new KPI, $\LBpace$, will guide both the learning process described in Section~\ref{sec:learning} and the illustration of the results in Section~\ref{sec:results}.

\subsection{Features} 
A key element of our approach is the identification and selection of input variables (features) that are able to explain and predict the behaviour of the aforementioned KPI for each branching rule. 
Following some of the ideas in \cite{Liberto2016}, \cite{Khalil2016}, and \cite{Alvarez2017} for MILP and our knowledge of \ref{eq:PO} problems, we have considered 34 features. They are summarized in Table \ref{table:features}. Features are categorized depending on whether they represent characteristics of the variables, constraints, monomials, coefficients, or other attributes of the polynomial optimization problems. Additionally, we have considered some novel features related to the graph representations of~\ref{eq:PO}, VIG and CMIG, as described in Section~\ref{sec:graph}. Many of the chosen features consist of counts or statistical measures (such as average, median, and variance) of the characteristics of the optimization problem. These are \textit{static} features, as they only depend on the formulation of the original problem. They represent global information of the problem instance, without any node-specific information that changes throughout the branch-and-bound tree.

\begin{table}[!htbp]
\caption{List of features used by the different learning techniques.}
\label{table:features}
\centering
\begin{threeparttable}[b]
\renewcommand{\arraystretch}{0.6}
{
    \begin{tabular}{ll}
      \toprule
      \multirow{5}{*}{Variables}   &  No. of variables, variance of the density of the variables\tnote{1} \\
         & Average/median/variance of the ranges of the variables\tnote{2}  \\
         & Average/variance of the no. of appearances of each variable\tnote{3} \\
         & Pct. of variables not present in any monomial with degree greater than one \\
         & Pct. of variables not present in any monomial with degree greater than two \\
        \midrule
      Constraints &  No. of constraints, Pct. of equality/linear/quadratic constraints\\
    \midrule
    \multirow{3}{*}{Monomials}   &  No. of monomials\\
      & Pct. of linear/quadratic monomials, Pct. of linear/quadratic RLT variables\\
      & Average pct. of monomials in each constraint and in the objective function\tnote{4}\\
      \midrule
      Coefficients   &  Average/variance of the coefficients\\
     \midrule
     \multirow{3}{*}{Other}   &  Degree and density of \ref{eq:PO}\tnote{5}\\
     & No. of variables divided by no. of constrains/degree \\
     & No. of RLT variables/monomials divided by no. of constrains\\
     \midrule
    Graphs  & Density, modularity, treewidth, and transitivity of VIG and CMIG\tnote{6}\\
    \bottomrule
    \end{tabular}}
\begin{tablenotes}
{\scriptsize
\item [1] The density of a variable is the number of different monomials in which it appears divided by the number of different monomials in the problem.
\item [2] The range of a variable is the difference between its upper and its lower bound.
\item [3] For each variable we compute in how many constraints (and objective function) that variable is present.
\item [4] For each constraint (and objective function) we compute the percentage of monomials present in it with respect to the total number of monomials of \ref{eq:PO}.
\item [5] The density of \ref{eq:PO} is its number of different monomials divided by total number of monomials a problem with the same degree an number of variables might have.
\item[6] We do not include transitivity for CMIG since, given the structure of the graph, its value is always 0.
\item[]
\vspace{-0.45cm}
}
\end{tablenotes}
\end{threeparttable}
\end{table}

\subsection{Quantile Regression Learning Techniques} \label{sec:learning}
In supervised learning, data consists of input–output pairs.
Given an optimization problem of the form \ref{eq:PO} we use, as input, the information from the features vector. Regarding the output, we compute for each of the branching rules described in Section~\ref{sec:branchingcriteria} the so-called $\LBpace$, that is, the ratio between the running time and the improvement in the lower bound. It is convenient to normalize the values of $\LBpace$ to $[0, 1]$ by dividing the best pace among all rules, \emph{i.e.}, the smallest, by the pace of each rule. We refer to this normalized pace as $\NLBpace$. Thus, the closer this value is to one for a given branching rule, the better the performance of that rule. Note that the $\LBpace$ of a given branching rule is independent of the other rules, but its $\NLBpace$ is not.

While our ultimate goal is the selection of the best branching rule, our variables of interest are quantitative outputs. Therefore, we address the problem as a regression problem. The task is to learn, for each rule, a real-valued function that models and predicts $\NLBpace$. Then, given an instance, the rule that corresponds to the highest predicted normalized pace is selected.

The primary goal of classical regression models is to estimate the conditional mean of a response variable, $y\in\mathbb{R}$, given the features vector $\mathbf{x}\in\mathbb{R}^d$. 
Let $y$ denote the $\NLBpace$ of a given rule. In the preliminary computational experiment in Figure~\ref{subfig:bp_comp} we observe an asymmetric behavior of $y$ (negative skewness), as well as the presence of outliers. This makes conventional regression models based on the conditional mean unsuitable. In this context, quantile regression is presented as a more appropriate methodology, since it does not make any assumptions about the distribution of the response variable and is more robust to the presence of outliers. Quantile regression aims to estimate  $Q_\tau(\mathbf{x})$, the $\tau$-conditional quantile of $y$, with $\tau\in (0,1)$. More precisely, 
$Q_\tau(\mathbf{x})=\inf\{y:\ F(y|\mathbf{x})\geq \tau\},$
being $ F(y|\mathbf{x})$ the conditional cumulative distribution function of $y$. Thus, the focus is on the conditional distribution of the response variable rather than just on its conditional mean, $\mu(\mathbf{x})=E(y|\mathbf{x})$. In our context it might be relevant to find out, for example, that for certain instances a given rule exceeds, with high probability, a given value of $y$. We refer to  \cite{Koenker2005} and \cite{Koenker2017}  for a detailed overview of quantile regression. 

There are several quantile regression methods in the statistics and ML literature. We focus on nonparametric methods, which allow for more flexible specifications of the relation between the response and the explanatory variables than their parametric counterparts. We present below a quick overview of the three main methods used in our analysis.


\subsubsection*{Quantile generalized additive models.}

Generalized additive models assume that the conditional mean of the response has an additive structure such as $\mu(\mathbf{x})=\sum_{j=1}^m f_j(\mathbf{x})$, where the $m$ additive terms represent predictors. The $f_j$'s are nonparametric functions which can be modeled using an expansion of basis functions; spline basis functions in our case. Hence, we get a flexible and interpretable statistical method, capable of characterizing non-linear regression effects. This leaves behind the traditional linear model that fails to represent the real effects of the inputs, which are normally non-linear. For detailed information of generalized additive models, we refer to  \cite{Hastie2009}. Quantile generalized additive models, see \cite{Fasiolo2021}, extend the previous idea and assume that the $\tau$-conditional quantile of $y$ is given as a linear combination of unknown smooth functions of the features vector. Thus, it uses a very flexible technique that combines the ideas of quantile regression and generalized additive models. 

\subsubsection*{Stochastic gradient boosting for quantile regression.}

Boosting is a supervised learning technique designed for classification problems, but also applicable to regression. The basic idea of this ensemble method is to create a highly accurate prediction rule as a weighted combination of weak learners predictions.
In this work we focus on the stochastic gradient boosting algorithm, introduced in \cite{Friedman2002}. The stochastic gradient boosting algorithm is a modification of the  gradient boosting algorithm by \cite{Friedman2001}, that improves the performance by incorporating randomness in the procedure. More specifically, a subsample of the training data is randomly selected at each iteration, and used to fit the base learner (tree).  
We use stochastic gradient boosting for quantile regression, which generalizes the described method to the context of quantile regression by using an asymmetrically
weighted absolute loss function \citep{Kriegler2010}. 

\subsubsection*{Quantile regression forests.}

Random forests are supervised learning algorithms that can be used both in classification and regression problems. They were introduced in \cite{Breiman2001} as ensemble methods that grow several individual decision trees and aggregate them to make a single prediction. For regression, they give an approximation of the conditional mean of a response variable. Quantile regression forests, introduced in \cite{Meinshausen2006}, are a generalization of random forests that compute an estimation of the conditional distribution of the response variable by taking into account all the observations in every leaf of every tree and not just their average. One of the advantages of random forests, as well as quantile random forests, is that we can use the complete dataset, without randomly splitting the data into training and test set, to evaluate the model using the Out-Of-Bag predictions.

\section{Computational Results}\label{sec:results}

\subsection{Experimental Setup} \label{sec:dataset} 
We use three different sets of instances to test the performance of the learning approach on a variety of instances with different characteristics. The first one is taken from \cite{Dalkiran2016} and consists of 180 instances of randomly generated polynomial optimization problems of different degrees, number of variables, and density. The second dataset comes from the well known benchmark MINLPLib \citep{minlplib}, a library of Mixed-Integer non-linear Programming problems. We have selected from MINLPLib those instances that are PO problems with box-constrained and continuous variables, resulting in a total of 168 instances. The third dataset comes from another well known benchmark, QPLIB \citep{qplib}, a library of quadratic programming instances, for which we made a selection analogous to the one made for MINLPLib, resulting in a total of 63 instances. Hereafter we refer to the first dataset as DS, to the second one as MINLPLib, and to the third one as QPLIB. Moreover, we have disregarded instances solved at the root node, since the branching rule plays no role in them. We also discarded unconstrained problems, as some of the features have the number of constraints in the denominator. Our final dataset consists of a total of 354 instances. 

All the executions reported in this paper have been performed on the supercomputer Finisterrae~II, at Galicia Supercomputing Centre (CESGA). Specifically, we used computational nodes powered with 2 deca-core Intel Haswell 2680v3 CPUs with 128GB of RAM and 1TB of hard drive. 

The branching rules defined in Section~\ref{sec:branchingcriteria} have been coded in RAPOSa \citep{Gonzalez-rodriguez2020}, a state-of-the-art implementation of the RLT outlined in Section~\ref{sec:PO}. The analyses are performed using a configuration of RAPOSa with a minimal number of enhancements (calls to an NLP local solver and J-sets \citep{Sherali2013}; refer to \cite{Gonzalez-rodriguez2020} for details). Each instance of our final dataset was run with six different configurations of RAPOSa, based on the six different branching rules defined in Section \ref{sec:branchingcriteria}. The time limit of each execution was set to one hour.

We  follow the standard learning procedure. We randomly split the complete set of instances into two disjoint sets: the training set (70\%) and the test set (30\%). With the objective of obtaining a better performance, we gather the instances into ``families'' related to the groupings defined in the corresponding libraries, where those belonging to the same group share similar characteristics. The within-family proportion of instances is maintained through the splitting process. Moreover, in order to gain robustness, we construct 10 partitions of the dataset into training and test data and report aggregate results over all the partitions.

As discussed in Section~\ref{sec:learning}, when performing quantile regression one has to specify the $\tau$-conditional quantile to be estimated. Due to the negative asymmetry present in the KPI (see Figure~\ref{fig:branchingcriteria}), we have opted for values below the median ($\tau\leq 0.5$), where more differentiated behaviors can be observed between the branching rules and where it is more relevant to make a good selection. Although the different values evaluated lead to similar results, the best ones were obtained with  $\tau=0.3$, which is therefore the value used for the numerical analysis.

\subsection{Numerical Results and Analysis}\label{subsec:results} 
This section presents the numerical results that allow to evaluate the performance of the obtained ML-based branching rules. The objective is to show the impact of the ML approach when embedded within RAPOSa and also to provide further evidence of its effectiveness across various sets of instances. As a first step we compare the three ML-based rules obtained with the learning methodologies described in Section~\ref{sec:learning}, that is, quantile generalized additive models (Q-GAM), stochastic gradient boosting for quantile regression (SGB-QR), and quantile regression forests (Q-RF). As RAPOSa includes several enhancements such as warm start and bound tightening \citep{Gonzalez-rodriguez2020}, we provide results of the ML-based branching without any of these enhancements (Baseline) and including them (WS and WS + BT). 

\begin{table}[!htbp]
    \caption{ML-based rules' performance with respect to $\LBpace$ (average over test sets).}
    \renewcommand{\arraystretch}{0.75}
    \centering
    \begin{tabular}{l|cccc}
    & \textbf{Baseline} & \textbf{WS} & \textbf{WS + BT} \\
    \toprule\toprule
    \textbf{Best branching rule} (across all instances) & 1.875 & 1.534 & 0.819 \\
    \textbf{ML-based branching rule:} & & & \\
    \multicolumn{1}{r|}{Q-GAM} & 1.439 & 1.181 & 0.782 \\
    \multicolumn{1}{r|}{SGB-QR} & 1.439 & 1.174 & 0.763 \\
    \multicolumn{1}{r|}{\textbf{Q-RF}}& 1.415  & 1.153 & 0.745 \\
    \textbf{Optimal rule} (instance by instance) & 1.187 & 0.956 & 0.653 \\
    \midrule
    \textbf{Improvement after learning} (with Q-RF) & 24.5\% & 24.8\% & 9.0\% \\
    \textbf{Optimal improvement} (upper bound for learning) & 36.7\% & 37.6\% & 20.2\% \\
     \bottomrule
    \end{tabular}
    \label{table:resultstestsets}
\end{table}

 Table~\ref{table:resultstestsets} presents the results for the $\LBpace$,\footnote{We take $\varepsilon=0.001$ in Equation~\eqref{eq:epsilon}.} where the reported number corresponds to the average, across the 10 test sets, of the geometric mean of $\LBpace$ in the instances of each test set. The best branching rule is the original branching rule where the aforementioned average is smaller, that is, the one performing best. For the Baseline configuration, the \texttt{range} is the best branching rule while for the other two setups, the \texttt{dual} is the best branching rule. The optimal rule corresponds with choosing, for each instance, the original rule that performs best. We can see that the three learning approaches Q-GAM, SGB-QR, and Q-RF have similar performance. Since Q-RF performs slightly better, hereafter we use it as the reference ML-based rule. For the Baseline configuration, compared to selecting the best branching rule, Q-RF improves the pace by 24.5\%. Importantly, this is not far from the improvement obtained by the optimal rule, 36.7\%, which is an upper bound for our algorithm selection approach. For the WS + BT, the improvement with Q-RF is 9\% and the optimal one is 20.2\%. This is due to the fact that bound tightening techniques already improve the bounds significantly, reducing the search space and, therefore, reducing as well the room for improvement using different branching rules. Interestingly, the WS enhancement improves the KPI from 1.875 to 1.534, which is 18.2\%. This value is lower than using Q-RF, which implies that given the choice between applying warm starting or a machine-learning branching approach, the latter may be better (and requiring no computational overhead). Table~\ref{table:resultsoob} shows that, when using out-of-bag predictions for the Q-RF methodology, we get similar results (the reported values correspond to the geometric means of $\LBpace$ across all the instances of the dataset and the best original rule was always the \texttt{dual}). 
\begin{table}[!htbp]
    \caption{ML-based rules' performance with respect to $\LBpace$ (Out-of-bag).}
\renewcommand{\arraystretch}{0.75}
    \centering
    \begin{tabular}{l|ccc}
    & \textbf{Baseline} & \textbf{WS} & \textbf{WS + BT} \\
    \toprule\toprule
    \textbf{Best branching rule} (across all instances) & 1.736 & 1.399 & 0.733 \\
    \textbf{Out-of-bag Q-RF} & 1.191 & 1.050 & 0.678 \\
    \textbf{Optimal rule} (instance by instance) & 0.983 & 0.810 & 0.601 \\
    \midrule
    \textbf{Improvement after learning} (with Q-RF) & 31.4\% & 24.9\% & 7.5\% \\
    \textbf{Optimal improvement} (upper bound for learning) & 43.4\% & 42.1\% & 18.0\% \\
     \bottomrule
    \end{tabular}
    \label{table:resultsoob}
\end{table}

\begin{figure}[!htbp]
\centering
\includegraphics[width=0.85\textwidth]{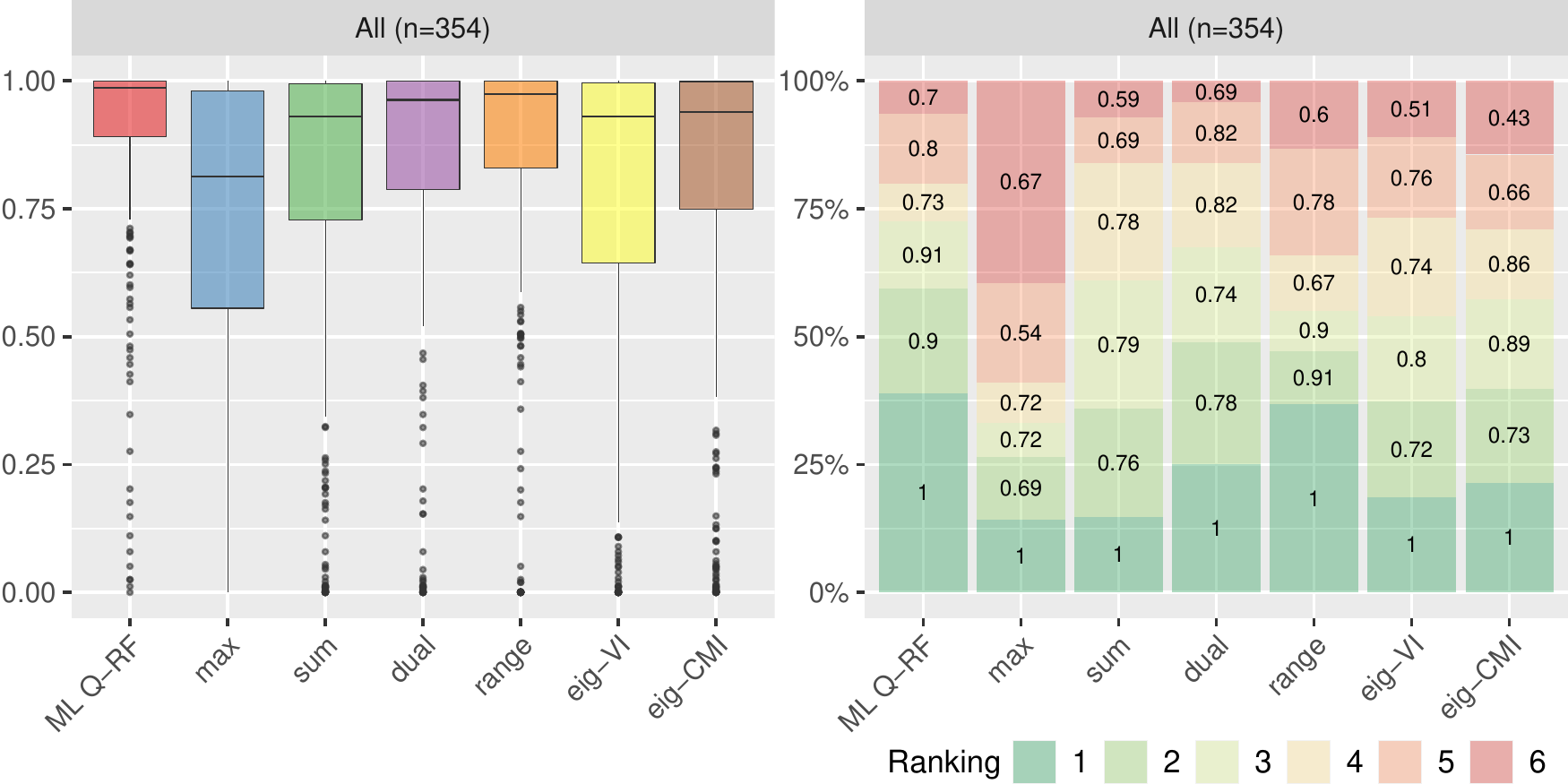}
\caption{Performance in Baseline configuration. Left, $\NLBpace$ boxplots of all the branching rules. Right, branching rules ranked from 1 (best) to 6 (worst), according to $\NLBpace$.}
\label{fig:noBTnoWS1h_results}
\end{figure}

In order to get additional insights from the above results, we now explore in depth the performance of the ML-based branching rule Q-RF in the Baseline setup, building upon the out-of-bag predictions.\footnote{We have also run the corresponding analysis for WS and WS+BT configurations, obtaining qualitatively similar results.} In Figure~\ref{fig:noBTnoWS1h_results}, the complete set of instances for the three libraries is used (total of 354 instances). The boxplots show that $\NLBpace$ is significantly closer to 1 for Q-RF than for any of the original rules. For the bar chart, it also validates that Q-RF is outperforming the other branching rules as it is ranked number 1 (best) in terms of the pace more often (38.98\% compared to range which is 36.72\%). The numbers inside the bars indicate average value of $\NLBpace$ for the problems in each rank position. Thus, even when Q-RF is ranked 6th (worst) its average normalized pace value is 0.7 which is also better than the corresponding value for all the other six branching rules (ranging from 0.43-0.69). In particular, the plots in Figure~\ref{fig:noBTnoWS1h_results} show that, although the learned criterion chooses the best rule in less than 50\% of the instances, it consistently chooses rules whose performance is very close to that of the best rule. This last observation also allows to understand why the improvement after learning with Q-RF, reported in Tables~\ref{table:resultstestsets} and~\ref{table:resultsoob}, is close to the optimal improvement.

\begin{figure}[!htbp]
    \centering
    \begin{subfigure}{\textwidth}
		\centering
        \includegraphics[width=0.85\textwidth]{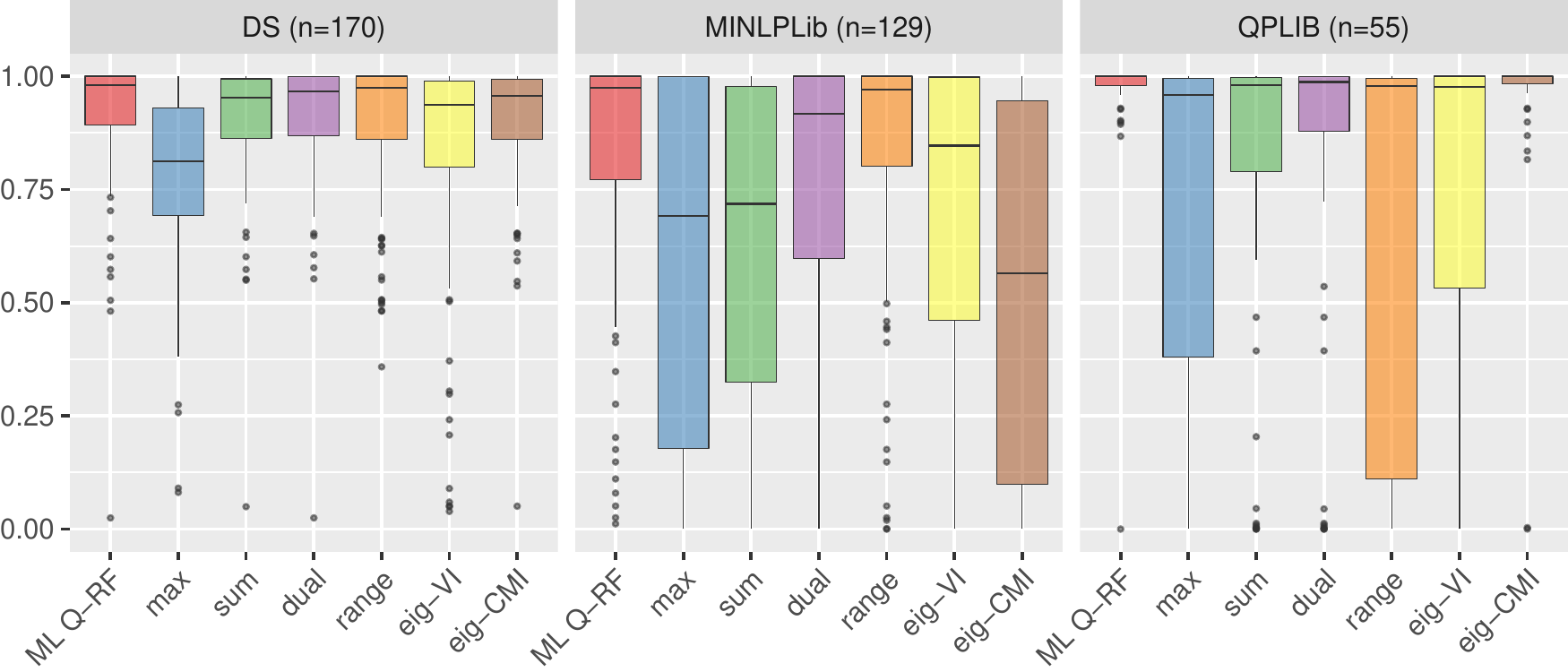}
        \caption{Boxplots of the normalized pace using the described branching rules.}
        \label{subfig:noBTnoWS1h_bp_ML}
    \end{subfigure}    
    
    \medskip
    \begin{subfigure}{\textwidth}
		\centering
        \includegraphics[width=0.85\textwidth]{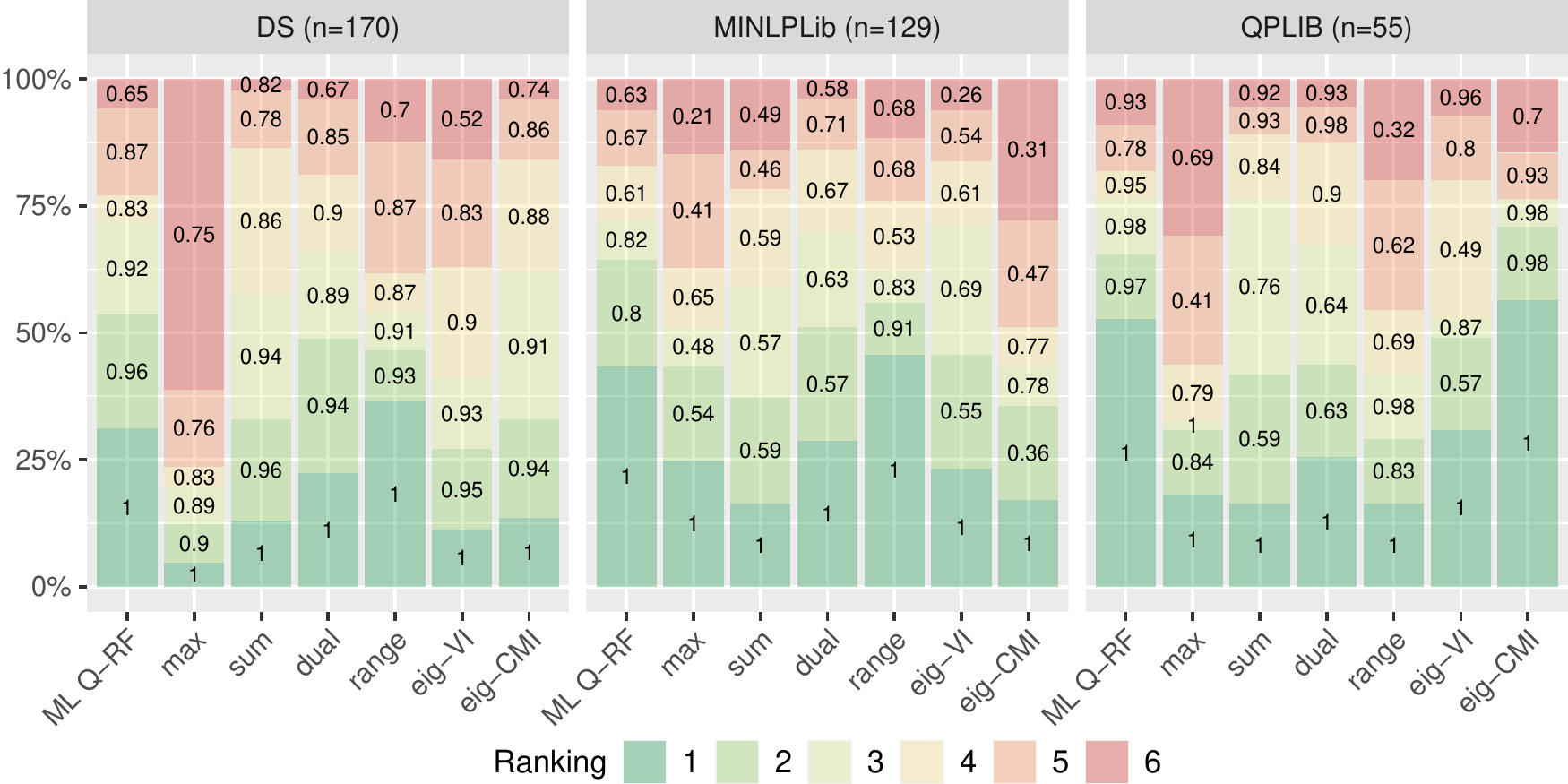}
        \caption{Branching rules ranked from 1 (best) to 6 (worst), according to the normalized pace.}
    \label{subfig:noBTnoWS1h_bar_ML}
    \end{subfigure}
    \caption{Performance in Baseline configuration. Comparison across the sets of instances.}
    \label{fig:noBTnoWS1h_byTS}
\end{figure}

In Figure~\ref{fig:noBTnoWS1h_byTS} we separate the results in Figure~\ref{fig:noBTnoWS1h_results} for the three sets of instances, DS, MINLPLib, and QPLIB. The number of instances within each library is indicated in parenthesis. There are some differences across datasets. For instance, the box plots show little variation in $\NLBpace$  across the branching rules for DS set compared to the other two sets. The \texttt{max} rule is the worst one, while the other rules are close in terms of normalized pace. For MINLPLib, \texttt{range} performs slightly better than Q-RF as the number of times that \texttt{range} ranked first is 45.73\% compared to 43.41\% for Q-RF. For QPLIB, both Q-RF and \texttt{eig-CMI} strongly dominate in terms of performance, ranking first 52.72\% and 56.36\% of the instances, respectively. Note that for this set, the \texttt{range} rule, which performed very well for MINLPLib, is the worst branching rule for QPLIB. In particular, it is ranked 6th in around 20\% of the instances with an average $\NLBpace$ of 0.32 in them, while the other rules, when ranked 6th, have an average $\NLBpace$ ranging from 0.69 to 0.96. Therefore, the ML-based rule seems to understand which are the features within and across datasets that drive the performance of the original rules.

\begin{figure}[!htbp]
\centering
\includegraphics[width=0.85\textwidth]{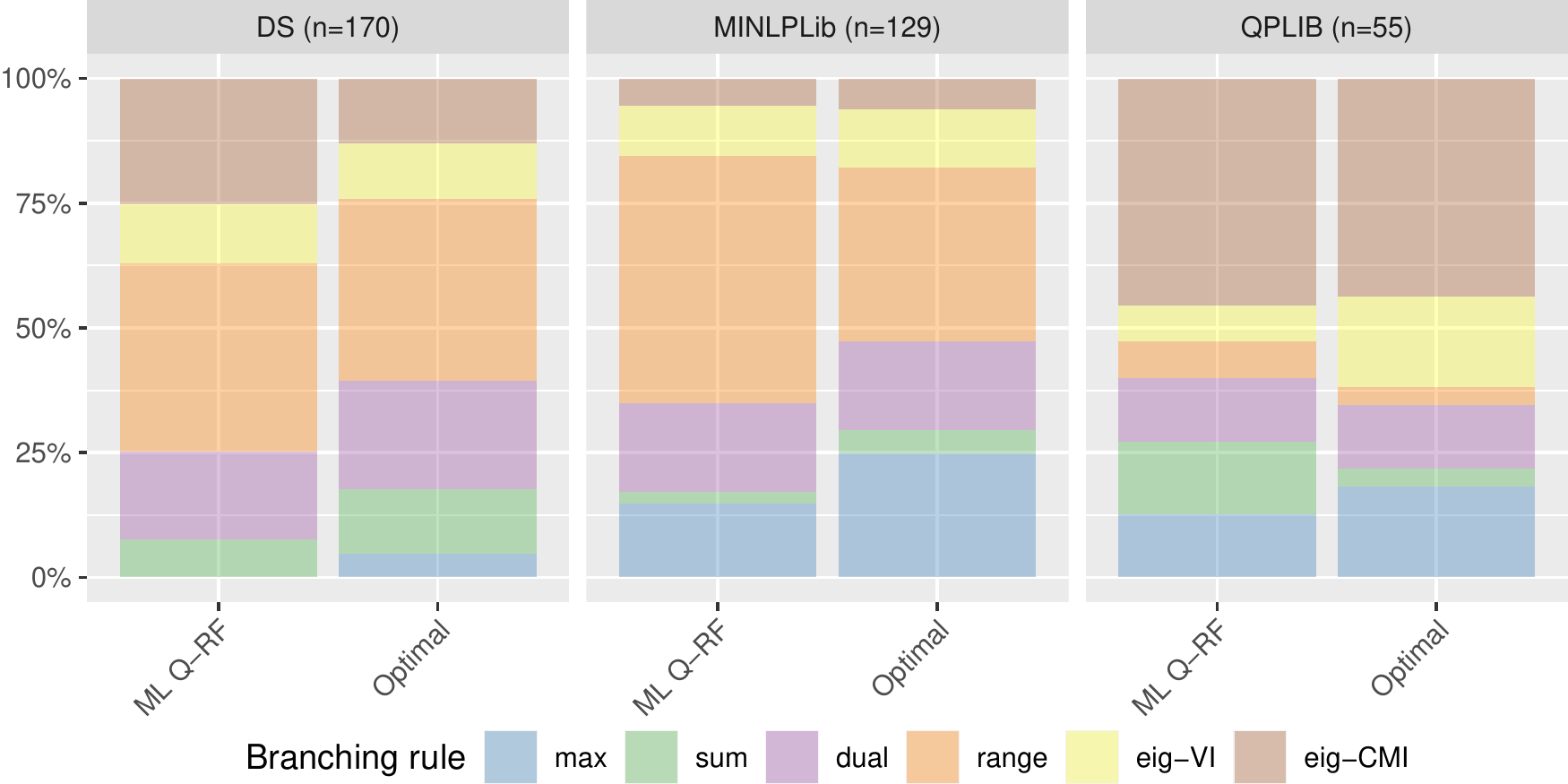}
\caption{Stacked bar graphs represent percentage of times that each branching rule is selected by Q-RF and the optimal branching rule.}
\label{fig:bar_comparison}
\end{figure}

Next, we further explore to what extent the ML-based rule manages to learn the patterns behind optimal rule selection. Figure \ref{fig:bar_comparison} presents a comparison of the frequencies with which the original rules are chosen by the ML-based rule and the optimal rule. The behavior of the ML-based branching rule strongly resembles that of the theoretical optimal branching rule in all three datasets. For DS the dominant rule is \texttt{range}, with \texttt{dual} and \texttt{eig-CMI} also playing an important role. For MINLPLib, \texttt{range} dominates again and, surprisingly, \texttt{max} rule is chosen quite often and so is again \texttt{dual}. For QPLIB, \texttt{eig-CMI} is chosen for almost 50\% of the instances for both the ML-based rule and the optimal rule. Interestingly, in this last set of instances, the \texttt{range} rule is rarely selected as the criteria for branching. Overall, these findings are consistent with the results presented in Figure~\ref{fig:noBTnoWS1h_byTS}. Recall that our offline learning is performed jointly on all sets of instances, so the results in Figure~\ref{fig:bar_comparison}, showing that the ML-based rule correctly adapts to the instances in DS, MINLPLib, and QPLib, are a sign of the quality of the learning process and of its potential to learn from highly hetereogeneous instances.

\begin{figure}[!htbp]
\centering
\includegraphics[width=0.85\textwidth]{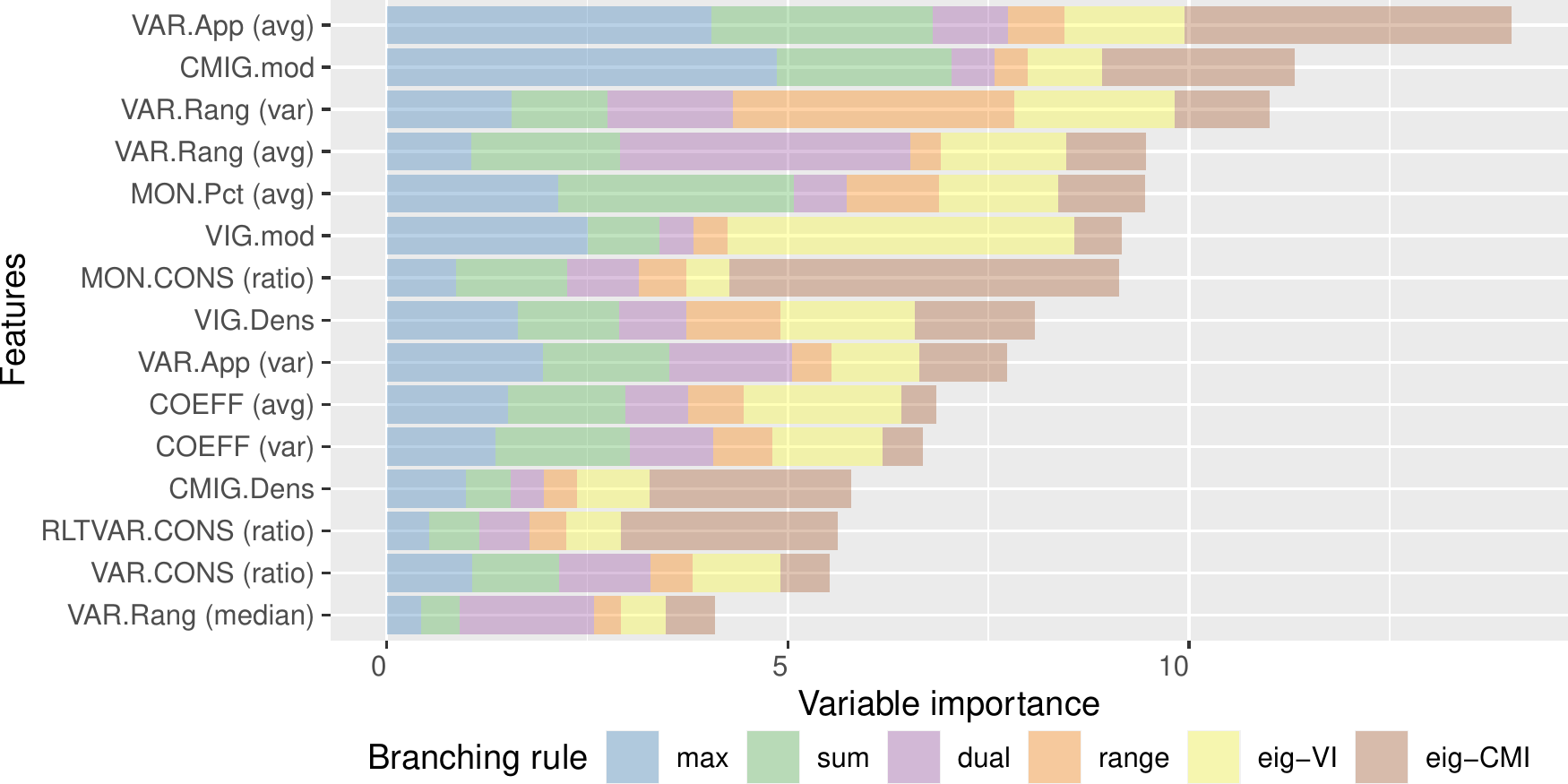}
\caption{Stacked bar graphs represent, for each branching rule, variable importance score of the features obtained using Q-RF.}
\label{fig:var_imp}
\end{figure}

One advantage of the learning techniques used in our analysis is the interpretability. In particular, (quantile) random forests allow to assign scores to the different features and, hence, understand the ones that are the most important for the associated predictions. Since our Q-RF builds upon six independent regressions, one for each of the original branching rules, we have six different scores for each explanatory feature. Figure~\ref{fig:var_imp} presents, for each of the six regressions, the feature importance scores of the top 15 features overall. A first observation is that the novel graph-related features are important, as three of them appear in the top 8: CMIG.mod, VIG.mod, and CMIG.Dens. Further, CMIG.mod, the modularity on CMIG, is in second position. Different features have different importance score for learning different branching rules. In particular, graph-related features are important for the corresponding graph-based rules: CMIG.mod and CMIG.Dens are particularly important for \texttt{eig-CMI} and VIG.mod  score is very high for \texttt{eig-VI}. Additionally, the most important features for learning \texttt{range} and \texttt{dual} branching rules are those related to the ranges of variables.

\begin{figure}[!htbp]
\centering
\includegraphics[width=0.85\textwidth]{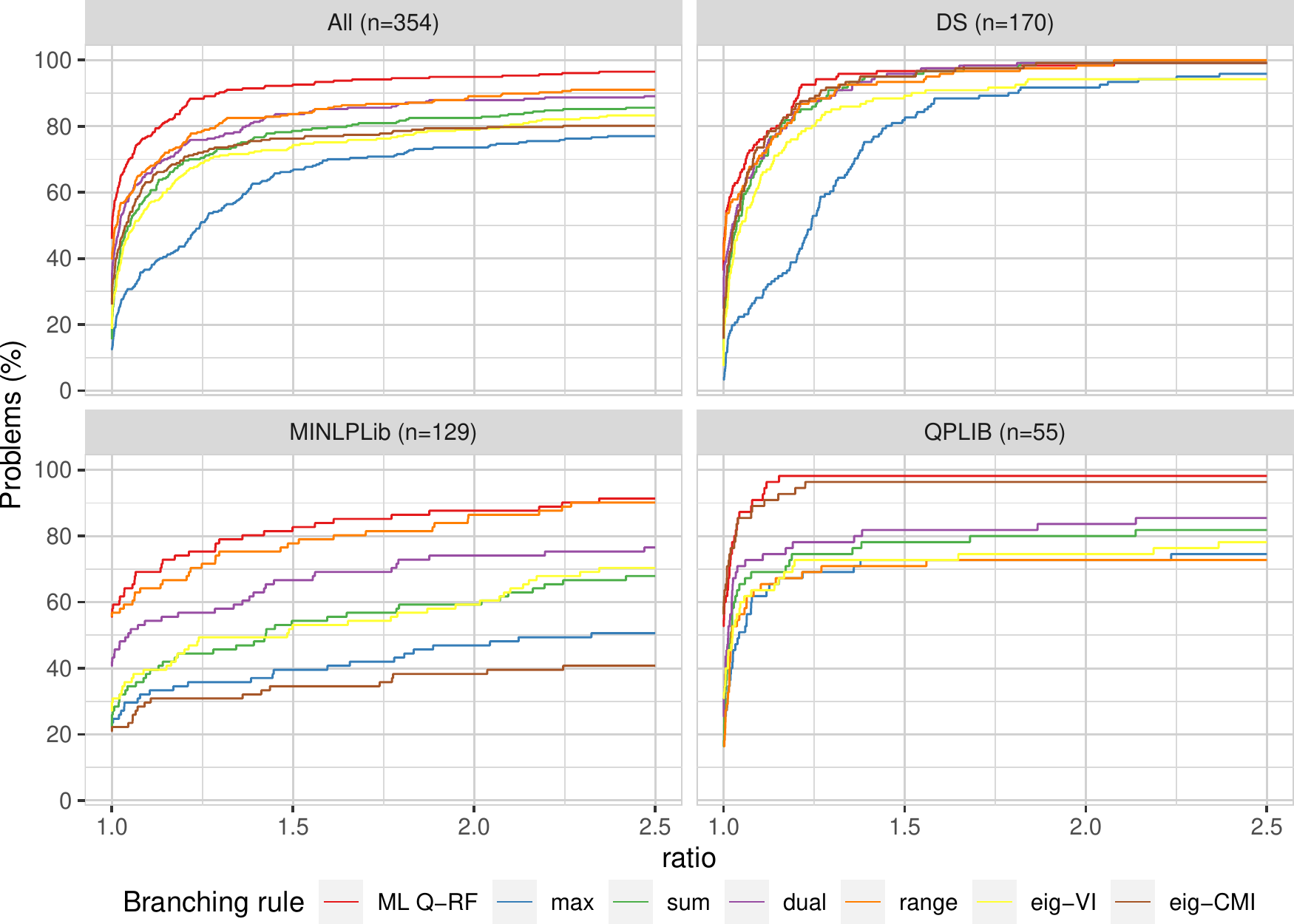}
\caption{Performance profiles using the pace KPI for different sets of instances.}
\label{fig:ppbat_ML}
\end{figure}

To conclude, we compare the performance of the different branching rules by showing the corresponding performance profiles \citep{Dolan2002}. In Figure~\ref{fig:ppbat_ML} we present the performance profiles for the three sets of instances together and individually. Consistently with all the results we have already reported so far, the ML-based rule Q-RF performs noticeably better than all the original branching rules. From the DS library plot, Q-RF performs slightly better than 4 of the 6 branching rules (while \texttt{max} and \texttt{eig-VI} are significantly worse). For the MINLPLib library, Q-RF is clearly the best and \texttt{range} is in second position. For the QPLIB, Q-RF is slightly better than \texttt{eig-CMI} and dominates the other 5 branching rules. Therefore, despite having that the original rules performing best vary significantly for the different sets of instances, Q-RF manages to slightly outperform the best rules in each individual set and solidly outperforms all the original rules overall.

\section{Conclusions and Future Research}\label{sec:conclusions} 

We have presented a framework for learning rules for spatial branching and shown in different numerical experiments that the resulting gains in performance can be significant (up to 30\%). These experiments have been performed jointly on a rich set of diverse instances taken from standard NLP libraries and from past literature. Interestingly, we have also shown that the interpretability of our framework allows to understand the importance of the different explanatory features for the resulting ML-based rule. This kind of understanding, which may shed light for the feature-design of future contributions, has been overlooked by most of the past literature on ``learning to branch''.

As we briefly discussed in Section~\ref{sec:framework}, the choice to perform learning only with static features limits the potential of our learning setting. Arguably, neglecting node by node information and limiting the scope of the algorithm selection to be instance specific, instead of node specific, may result in learning a less effective branching rule. Yet, given that this paper is the first contribution on ``learning to branch'' in the non-linear setting, we believe that the static learning setting provides a natural starting point. First, since it involves no computational overhead for the computation of node-specific features, the differences in performance of the learned rule with respect to the original ones directly reflects the net effect of the learning process. Second, learning static rules allows for a clean definition of the ``optimal'' rule for each instance (this would not be entirely straightforward with dynamic rules, where the selection of the best ``myopic'' rule in each node might not lead to the best performance on the given instance). This is relevant since, once we pin down what can be achieved with perfect learning, one can easily assess the quality of the learning behind the rules resulting from our approach. Third, static learning provides a benchmark on which to assess the performance of more advanced dynamic techniques and, in particular, to evaluate to what extent the online computational overhead of node-specific features is compensated by the reduction on the tree size.

Given the above discussion, a prominent direction for future research is the incorporation of dynamic features and study the performance of the resulting rules taking the approach developed in this paper as a benchmark. Additionally, future research should aim at having richer portfolios of branching rules, incorporating for instance rules that build upon reliability branching or violation transfer (see, for instance, \cite{Tawarmalani2004,Achterberg2005,Belotti2009}).  

{\footnotesize
\ACKNOWLEDGMENT{
This research has been funded by FEDER and the Spanish Ministry of Science and Technology through projects MTM2014-60191-JIN and MTM2017-87197-C3. Brais Gonz\'alez-Rodr\'iguez and Ignacio G\'omez-Casares acknowledge support from the Spanish Ministry of Education through FPU grants 17/02643 and 20/01555, respectively. Beatriz Pateiro-L\'{o}pez acknowledges support from Grant PID2020-116587GB-I00 funded by MCIN/AEI/10.13039/501100011033 and ED431C 2021/24 funded by Conseller\'ia de Cultura, Educaci\'{o}n e Universidade.  Bissan Ghaddar's research is supported by Natural Sciences and Engineering Research Council of Canada Discovery Grant 2017-04185 and by the David G. Burgoyne Faculty Fellowship.}}

\bibliographystyle{apalike}
\addcontentsline{toc}{chapter}{\bibname}
\bibliography{references}

\end{document}